\documentclass[11pt]{article}

\usepackage{amssymb}
\usepackage{amsmath}

\def\square{\hbox{\vrule\vbox{\hrule\phantom{o}\hrule}\vrule}}
\def\endofproofsymbol{\square}

\makeatletter \def\@cite#1#2{{#1\if@tempswa , #2\fi}}
              \def\@biblabel#1{#1.}
\makeatother

\newtheorem{definition}{Definition}[section]

\begin{document}

\title{Contrasting Two Transformation-Based Methods
for Obtaining Absolute Extrema}

\author{\sc{D. F. M. Torres}$^1$
\sc{and} \sc{G. Leitmann}$^2$}

\date{}

\maketitle

%%%%%%%%%%%%%%%%%%%%%%%%%%%%%%%%%%%%%%%%%%%%

\footnotetext[1]{Associate Professor in the Department of Mathematics,
University of Aveiro, Aveiro, Portugal.}

\footnotetext[2]{Professor in the Graduate School,
College of Engineering, University of California,
Berkeley, California.}

%%%%%%%%%%%%%%%%%%%%%%%%%%%%%%%%%%%%%%%%%%%%

\paragraph*{Abstract.}
In this note we contrast two transformation-based methods to
deduce absolute extrema and the corresponding extremizers. Unlike
variation-based methods, the transformation-based ones of Carlson
and Leitmann and the recent one of Silva and Torres are direct in
that they permit obtaining solutions by inspection.

\bigskip


\paragraph*{Key~Words.} Optimal Control, direct methods,
Calculus of Variations, absolute extrema, invariance, symmetry.

%%%%%%%%%%%%%%%%%%%%%%%%%%%%%%%%%%%%%%%%%%%%

\section{Introduction}

In the mid 1960's a direct method for the problems of the calculus
of variations, which permits one to obtain absolute extremizers
directly, without using variational methods, was introduced by
Leitmann (Ref.~\cite{Leitmann67}). Since then, this direct method has
been extended and applied to a variety of problems (see \textrm{e.g.}
Refs.~\cite{Carlson02,CarlsonLeitmann05a,CarlsonLeitmann05b,Leitmann01}).
A different but related direct approach to problems of optimal
control, based on the variational symmetries of the problem
(\textrm{cf.} Refs.~\cite{GouveiaTorres05,GouveiaTorresRocha06}), was
recently introduced by Silva and Torres (Ref.~\cite{SilvaTorres06}).
The emphasis in Ref.~\cite{SilvaTorres06} has been on showing the
differences and similarities between the proposed method and that
suggested by Leitmann. In order to illustrate the relation between
these two methods, only examples capable of treatment by both
methods were presented in Ref.~\cite{SilvaTorres06}. In this note,
we discuss some differences between the method of Carlson and
Leitmann (C/L) and Silva and Torres (S/T). In particular, we show
how one succeeds when the other does not.

%%%%%%%%%%%%%%%%%%%%%%%%%%%%%%%%%%%%%%%%%%%%%%%%%%%%%%

\section{The Invariant Transformation Method of S/T$^3$}
\label{sec:ITM}

\footnotetext[3]{Throughout this Note, the notation conforms
to that used in the references.}

Let us consider the problem of optimal control in Lagrange form:
minimize an integral
\begin{equation}
\label{probCO1} I\left[x(\cdot), u(\cdot)\right] = \int_a^b L
\left(t, x(t), u(t) \right) dt
\end{equation}
subject to a control system
\begin{equation}
\label{controlsystCO1} \dot{x}(t) = \varphi \left(t, x(t),
u(t)\right) \quad \text{a.e. on } [a, b] \, ,
\end{equation}
together with appropriate boundary conditions and constraint on
the values of the control variables:
\begin{equation}
\label{eq:const} x(a) = x_a \, , \quad x(b) = x_b \, , \quad u(t)
\in \Omega \, .
\end{equation}
The Lagrangian $L(\cdot, \cdot, \cdot)$ is a real function assumed
to be continuously differentiable in $[a,b] \times \mathbb{R}^{n}
\times \mathbb{R}^{m}$; $t \in \mathbb{R}$ is the independent
variable; $x(\cdot) : [a,b] \rightarrow \mathbb{R}^{n}$ the vector
of state variables; $u(\cdot) : [a,b] \rightarrow \Omega \subseteq
\mathbb{R}^{m}$ the vector of controls, assumed to be a piecewise
continuous function; and $\varphi(\cdot,\cdot,\cdot) : [a,b]
\times \mathbb{R}^{n} \times \mathbb{R}^{m} \rightarrow
\mathbb{R}^{n}$ the derivative function, assumed to be a
continuously differentiable vector function. When $\Omega$ is an
open set (it may be all Euclidean space $\Omega = \mathbb{R}^m$),
problem \eqref{probCO1}-\eqref{eq:const} can be studied using the
classical techniques of the Calculus of Variations. Optimal
Control Theory includes the classical Calculus of Variations and
generalizes the theory by dealing with the cases where $\Omega$ is
not an open set.

The application of the invariant transformation method
(Ref.~\cite{SilvaTorres06}) depends on the existence of a sufficiently
rich family of invariance transformations (variational
symmetries). The reader interested on the study of variational
symmetries is referred to Refs.~\cite{GouveiaTorresRocha06,Torres04,Torres06}
and references therein.

\begin{definition}
\label{definv} Let $h^s$ be a $s$-parameter family of $C^1$
mappings satisfying:
\begin{equation*}
\begin{split}
&h^s(\cdot,\cdot,\cdot) : [a,b]\times \mathbb{R}^{n} \times \Omega
\longrightarrow
\mathbb{R} \times \mathbb{R}^{n} \times \mathbb{R}^{m} \, , \\
& h^s(t,x,u) = \left( t^s(t,x,u), x^s(t,x,u),
u^s(t,x,u) \right) \, , \\
& h^0(t,x,u)= (t,x,u) \, , \quad \forall (t,x,u) \in [a,b] \times
\mathbb{R}^{n} \times \Omega \, .
\end{split}
\end{equation*}
If there exists a function $\Phi^s(t,x,u) \in C^1\left([a,b],
\mathbb{R}^{n}, \Omega; \mathbb{R} \right)$ such that
\begin{equation}
\label{eq:inv:L} L \circ h^s(t, x(t), u(t)) \frac{d}{dt}t^s
\left(t, x(t), u(t)\right) = L\left(t, x(t), u(t)\right) +
\frac{d}{dt}\Phi^s\left(t, x(t), u(t)\right)
\end{equation}
and
\begin{equation}
\label{eq:inv:CS} \frac{d}{dt}x^s\left(t, x(t), u(t)\right) =
\varphi \circ h^s\left(t, x(t), u(t) \right) \frac{d}{dt}
t^s\left(t, x(t), u(t)\right)
\end{equation}
for all admissible pairs $\left(x(\cdot), u(\cdot)\right)$, then
\eqref{probCO1}--\eqref{controlsystCO1} is said to be
\emph{invariant} under the transformations $h^s(t,x,u)$ up to
$\Phi^s(t,x,u)$; and the transformations $h^s(t,x,u)$ are said to be
a \emph{variational symmetry} of \eqref{probCO1}--\eqref{controlsystCO1}.
\end{definition}

The method proposed in Ref.~\cite{SilvaTorres06} is based on a very
simple idea. Given an optimal control problem, one begins by
determining its invariance transformations according to
Definition~\ref{definv}. With respect to this, the tools developed
in Refs.~\cite{GouveiaTorres05,GouveiaTorresRocha06} are
useful. Applying the parameter-invariance transformations, we embed our
problem into a parameter-family of optimal control problems. Given
the invariance properties, if we are able to solve one of the
problems of this family, we also get the solution to our original
problem (or to any other problem of the same family) from the
invariant transformations. In section~\ref{sec:4} we give an
example which shows that the Invariant Transformation Method
(Ref.~\cite{SilvaTorres06}) is more general than the earlier C/L
transformation method in the case of optimal control problems.

%%%%%%%%%%%%%%%%%%%%%%%%%%%%%%%%%%%%%%%%%%%%%%%%%%%%%%

\section{The Direct Solution Method of C/L}

Since this method is fully discussed in readily available references, \textrm{e.g.}
Refs.~\cite{Carlson02,CarlsonLeitmann05a,CarlsonLeitmann05b,Leitmann67,Leitmann01},
many in this journal, we shall only recall that
the C/L transformation based method is applicable to problems
in the Calculus of Variations format: minimize an integral
\begin{equation}
\label{eq:6} I\left[x(\cdot)\right] = \int_a^b F \left(t, x(t),
\dot{x}(t) \right) dt
\end{equation}
with given end conditions
\begin{equation}
\label{eq:7} x(a) = x_a \, , \quad x(b) = x_b \, .
\end{equation}
If one wishes to solve an optimal control problem
\eqref{probCO1}-\eqref{eq:const}, the ``elimination'' of $u(t)$ in
favor of a function of $t$, $x(t)$, $\dot{x}(t)$ must be possible.
As illustrated in section~\ref{sec:4}, this may fail even if the
Implicit Function Theorem is satisfied.

Both the S/T and the C/L methods are predicated on posing a
problem ``equivalent'' to the original problems
\eqref{probCO1}-\eqref{eq:const} and \eqref{eq:6}-\eqref{eq:7},
respectively. Thus, these methods are useful only if the solution
of the ``equivalent'' problem is \emph{directly} obtainable,
\textrm{i.e.}, by inspection. There is, at present, no result
assuring that this can be done in general for the S/T method.
However, for the C/L method, at least in the scalar $x$ case, it
has been shown in Ref.~\cite{TE} and generalized to open-loop
differential games in Ref.~\cite{DE}, that the ``equivalent''
problem always has a minimizing solution obtained by inspection.
The conditions sufficient for this result are convexity of
integrand $F \left(t, x(t), \dot{x}(t) \right)$ with respect to
$\dot{x}(t)$, and existence of a so-called ``field of extremals''.
Indeed, no matter what the integrand of the original problem is,
provided the conditions above are met, the absolute minimizer of
the equivalent problem is always a constant.

%%%%%%%%%%%%%%%%%%%%%%%%%%%%%%%%%%%%%%%%%%%%%%%%%

\section{Example 1}
\label{sec:4}

The advantage of the invariant transformation method when compared
with the earlier transformation method is that one can apply it
directly to control systems whereas the method of C/L requires
that the control $u(t)$ can be expressed as a sufficiently smooth
function of $t$, $x(t)$, $\dot{x}(t)$, \textrm{e.g.} such that the
integrand be continuous in $x(t)$ and $\dot{x}(t)$. Here we use
the invariant transformation method of S/T to solve a simple
optimal control problem that is not covered by the classical
theory of the Calculus of Variations and which can not be solved
by the previous transformation method.

Consider the global minimum problem
\begin{equation}
\label{ex1}
\begin{gathered}
I[u_1(\cdot),u_2(\cdot)] = \int_0^1
\left(u_1(t)^2 + u_2(t)^2\right) dt \longrightarrow \min \\
\begin{split}
\dot{x}_1(t) &= \exp(u_1(t)) + u_1(t) + u_2(t) \, ,\\
\dot{x}_2(t) &= u_2(t) \, ,
\end{split}\\
x_1(0) = 0 \, , \quad x_1(1) = 2 \, , \quad
x_2(0) = 0 \, , \quad x_2(1) = 1 \, , \\
u_1(t) \, , u_2(t) \in \Omega = [-1,1] \, .
\end{gathered}
\end{equation}
We apply the procedure introduced in Ref.~\cite{SilvaTorres06} and
briefly described in section~\ref{sec:ITM}. First we notice that
problem \eqref{ex1} is variationally invariant according to
Definition~\ref{definv} under the one-parameter
transformations$^4$
\footnotetext[4]{A computer algebra package
that can be used to find the invariance transformations
(see Refs.~\cite{GouveiaTorres05,GouveiaTorresRocha06})
is available from the \emph{Maple Application Center} at
\texttt{http://www.maplesoft.com/applications/app\_center\_view.aspx?AID=1983}}
\begin{equation}
\label{transfEx1}
x_1^s = x_1 + s t \, , \quad x_2^s = x_2 + s t \, ,
\quad u_2^s = u_2 + s \quad (t^s = t \text{ and } u_1^s = u_1) \, .
\end{equation}
To prove this, we need to show that both the functional integral $I[\cdot]$
and the control system stay invariant under the $s$-parameter
transformations \eqref{transfEx1}. This is easily seen
by direct calculations. We begin by showing \eqref{eq:inv:L}:
\begin{equation}
\label{inv:Func:Ex1}
\begin{split}
I^s[u_1^s(\cdot),u_2^s(\cdot)]&=
\int_0^1 \left(u_1^{s}(t)\right)^2 + \left(u_2^{s}(t)\right)^2  dt \\
&= \int_0^1 u_1(t)^2 + \left(u_2(t) + s \right)^2  dt \\
&= \int_0^1 u_1(t)^2 + u_2(t)^2  + \frac{d}{dt} \left(s^2 t + 2 s x_2(t) \right) dt \\
&= I[u_1(\cdot),u_2(\cdot)] + s^2 + 2s \, .
\end{split}
\end{equation}
We remark that $\Phi^s\left(t,x_2\right) = s^2 t + 2s x_2$
and that $I^s[\cdot]$ and $I[\cdot]$ have the same minimizers:
adding a constant $s^2 + 2s$ to the functional $I[\cdot]$ does not
change the minimizer of $I[\cdot]$. It remains to prove \eqref{eq:inv:CS}:
\begin{equation}
\label{inv:CS:Ex1}
\begin{split}
\frac{d}{dt} \left(x_1^{s}(t)\right)&= \frac{d}{dt}
\left(x_1(t) + s t \right) = \dot{x}_1(t) + s = \exp(u_1(t)) + u_1(t) + u_2(t) + s \\
&= \exp(u_1^s(t)) + u_1^s(t) + u_2^s(t) \, , \\
\frac{d}{dt} \left(x_2^{s}(t)\right) &= \frac{d}{dt}
\left(x_2(t) + st \right) = \dot{x}_2(t) + s = u_2(t) + s \\
&= u_2^{s}(t) \, .
\end{split}
\end{equation}
Conditions \eqref{inv:Func:Ex1} and \eqref{inv:CS:Ex1} prove that
problem \eqref{ex1} is invariant under the $s$-parameter
transformations \eqref{transfEx1} up to the gauge term $\Phi^s =
s^2 t + 2s x_2$. Using the invariance transformations
\eqref{transfEx1}, we generalize problem \eqref{ex1} to a
$s$-parameter family of problems, $s \in \mathbb{R}$, which
include the original problem for $s=0$:
\begin{equation*}
\begin{gathered}
I^s[u_1^s(\cdot),u_2^s(\cdot)]
= \int_0^1 (u_1^s(t))^2 + (u_2^s(t))^2 dt \longrightarrow \min \\
\begin{split}
\dot{x}_1^s(t) &= \exp(u_1^s(t)) + u_1^s(t) + u_2^s(t) \, ,\\
\dot{x}_2(t) &= u_2^s(t) \, ,
\end{split}\\
x_1^s(0) = 0 \, , \quad x_1^s(1) = 2 + s \, , \quad
x_2^s(0) = 0 \, , \quad x_2^s(1) = 1 + s \, , \\
u_1^s(t) \in [-1,1] \, ,  \quad u_2^s(t) \in [-1+s,1+s] \, .
\end{gathered}
\end{equation*}
It is clear that $I^s \geq 0$ and that $I^s=0$ if
$u_1^{s}(t)=u_2^{s}(t) \equiv 0$. The control equation, the
boundary conditions and the constraints on the values of the
controls, imply that $u_1^{s}(t)=u_2^{s}(t) \equiv 0$ is
admissible only if $s=-1$: $x_1^{s=-1}(t) = t$, $x_2^{s=-1}(t)
\equiv 0$. Hence, for $s= -1$ the global minimum to $I^s[\cdot]$
is 0 and the minimizing trajectory is given by
\begin{equation*}
\tilde{u}_1^{s}(t) \equiv 0  \, , \quad \tilde{u}_2^{s}(t) \equiv
0 \, , \quad \tilde{x}_1^{s}(t)= t  \, , \quad \tilde{x}_2^{s}(t)
\equiv 0 \, .
\end{equation*}
Since for any $s$ one has by \eqref{inv:Func:Ex1} that
$I[u_1(\cdot),u_2(\cdot)] = I^s[u_1^s(\cdot),u_2^s(\cdot)] - s^2 - 2s$,
we conclude that the global minimum for problem $I[u_1(\cdot),u_2(\cdot)]$ is 1.
Thus, using the inverse functions of the variational symmetries
\eqref{transfEx1},
\begin{equation*}
u_1(t) = u_1^{s}(t) \, , \quad u_2(t) = u_2^{s}(t)- s \, , \quad
x_1(t) = x_1^{s}(t) - s t \, , \quad x_2(t) = x_2^{s}(t) - s t \, ,
\end{equation*}
the absolute minimizer is
\begin{equation*}
\tilde{u}_1(t) = 0 \, , \quad \tilde{u}_2(t) = 1 \, , \quad
\tilde{x}_1(t) = 2 t \, , \quad \tilde{x}_2(t) = t \, .
\end{equation*}
This problem cannot be solved by the C/L method.

While the existence of a function $h(\cdot)$ such that $u_1(t) =
h\left(\dot{x}_1(t)-\dot{x}_2(t)\right)$ is assured by the
satisfaction of the Implicit Function Theorem, this is not useful
in ``eliminating'' the control in favor of $t$, $x(t)$, $\dot{x}(t)$
which is a requirement of the C/L method. This is so because
$h(\cdot)$ is a solution of a transcendental equation. In
addition, since the controls are bounded, even if ``elimination''
of the control were possible, this would lead to differential
constraints of the form briefly discussed in Ref.~\cite{UE}.

%%%%%%%%%%%%%%%%%%%%%%%%%%%%%%%%%%%%%%%%%%%%%%

\section{Example 2}

Consider the problem of attaining the absolute minimum of integral
\begin{equation}
\label{eq:12}
I[x(\cdot)] = \int_a^b \left[ \dot{x}^2(t) + t \dot{x}(t)\right] dt
\end{equation}
with prescribed end conditions
\begin{equation}
\label{eq:13}
x(a) = x_a \, , \quad x(b) = x_b \, .
\end{equation}
This is a problem in the Calculus of Variations format for which
the C/L or S/T transformation-based method applies. These methods
can be applied \textit{ab initio} to obtain the solution. However,
here we shall employ the constructive sufficiency condition
embodied in Theorem~7 of Ref.~\cite{TE} towards that end.

The Euler-Lagrange equation is
\begin{equation}
\label{eq:14}
\ddot{x}(t) = -\frac{1}{2}
\end{equation}
so that
\begin{equation}
\label{eq:15}
x(t) = -\frac{1}{4} t^2 + c_1 t + c_2
\end{equation}
is the extremal with the constants of integration given by end
conditions \eqref{eq:13}, say $c_i = c_i^{*}$, $i=1, 2$, and the
extremal satisfying \eqref{eq:15} is
\begin{equation}
\label{eq:16}
x^{*}(t) = -\frac{1}{4} t^2 + c_1^{*} t + c_2^{*}
\end{equation}
which, being the solution of necessary condition \eqref{eq:14}, is
a candidate for the absolute minimizer of
\eqref{eq:12}-\eqref{eq:13}.

Now we can readily show that the conditions of Theorem~7 of
Ref.~\cite{TE} are met. Consider the one-parameter family of
extremals
\begin{equation*}
\xi(t,\beta) = -\frac{1}{4} t^2 + c_1^{*} t + c_2^{*} + \beta \, .
\end{equation*}
First of all, there exist a $\beta^{*}$, namely $\beta^{*} = 0$,
such that
\begin{equation*}
\xi(t,\beta^{*}) = x^{*}(t) \, .
\end{equation*}
Secondly,
\begin{equation*}
\frac{\partial \xi(t,\beta)}{\partial \beta} \ne 0 \, ,
\end{equation*}
and finally, the integrand of \eqref{eq:12} is convex in
$\dot{x}(t)$.

Thus, the extremal \eqref{eq:16} is indeed the absolute minimizer
of \eqref{eq:12}-\eqref{eq:13}.

Of course, this is precisely the solution obtained by employing
the C/L method. Indeed, the more general method inherent in
Theorem~7 of Ref.~\cite{TE} was derived using the C/L
transformation-based method.

%%%%%%%%%%%%%%%%%%%%%%%%%%%%%%%%%%%%%%%%%%%%%%

\section{Conclusion}

We have contrasted two transformation-based methods for obtaining
absolutely extremizing solutions for two classes of problems. One,
dubbed the Carlson/Leitmann method, is applicable to problems in
the format of the Calculus of Variations. The other, due to Silva
and Torres, is applicable to problems of Optimal Control.

We have shown that it is not always possible to convert an Optimal
Control problem into a Calculus of Variations one. Hence, the C/L
method may fail to apply when the S/T succeeds. On the other hand,
a classical constructive sufficiency condition, readily derived by
the C/L method, renders absolute extremizers for specific problems
of the Calculus of Variations more directly and succinctly than
the C/L and S/T methods.

%%%%%%%%%%%%%%%%%%%%%%%%%%%%%%%%%%%%%%%%%%%%%%

%%%%%%%%%%%%%%%%%%%%%%%%%%%%%%%%%%%%%%%%%%%%%%%%%%%%%%%

\end{document}